\documentclass[12pt]{article}
\usepackage[english]{babel}
\usepackage{amsmath,amsthm}
\usepackage{amsfonts}
\usepackage{amssymb}

\theoremstyle{definition}

\theoremstyle{remark}

\numberwithin{equation}{section}

\begin{document}

\begin{center}
\bfseries  A Formula of the Dirichlet Character Sum  \\
\end{center}

\begin{center}

\vspace{5mm}
JinHua  Fei\\
\vspace{5mm}

ChangLing  Company of Electronic Technology    \, Baoji \,  Shannxi  \,  P.R.China \\
\vspace{5mm}

E-mail: feijinhuayoujian@msn.com \\

\end{center}

\vspace{3mm}
{\bfseries Abstract.}\, In this paper, We use the Fourier series expansion of real variables function,  We give a formula to calculate the Dirichlet character sum, and four special examples are given.


{\bfseries MR(2000)  Subject  Classification \quad 11L40 } \\

\vspace{8mm}

The calculation of the Dirichlet character sum is very important in the number theory. This paper uses the Fourier series expansion of the functions, we give a general formula for calculating the Dirichlet character sum, Then, four examples are given to illustrate.

In this paper, $ \chi_q  $ denote the Dirichlet primitive character of mod $q $, If $ f(x) $ is a real function, we write

$$ f^*(x)=\frac{f(x+0)+f(x-0)}{2}  $$

and

$$  G(n,\chi_q) = \sum^{q-1}_{k=1} \chi_q (k) e\left(\frac{kn}{q}\right)  \qquad   \tau (\chi_q) = \sum^{q-1}_{k=1} \chi_q (k) e\left(\frac{k}{q}\right)   $$ \\

where $ e(x) = e^{2\pi i x} $ \\

First, let's give some lemmas.\\

{\bfseries Lemma 1.}   If $ \chi_q $ is the primitive character of module $q$,  then we have

$$G(n,\chi_q) = \sum^{q-1}_{k=1} \chi_q (k) e\left(\frac{kn}{q}\right)  = \overline{\chi}_q(n) \tau (\chi_q)  $$\\

see page 287 of references[1].\\

{\bfseries Lemma 2.}  If $ \chi_q $ is the primitive real character of module $q$,  then we have\\

\[   \tau(\chi_q )  = \left\{
  \begin{array}{ll}
  \sqrt{q} \quad if \, \, \, \chi_q (-1)=1 \\\\

  i\sqrt{q} \quad  if \, \, \, \chi_q (-1)=-1    \qquad \qquad \qquad
  \end{array}
  \right.
\]\\

see page 167 of references[2].\\

{\bfseries Lemma 3.}  If  the function  $ f(x) $ is defined in the interval $  [0,1]  $   and satisfies the  Dirichlet condition,  then we have

$$  f^*(x)=\int_0^1 f(t) dt +  \sum_{  n =-\infty \atop n  \neq 0 }^{\infty} \exp(2\pi i nx ) \int_0^1 f(t) \exp(- 2\pi i nt )  dt $$ \\

see page 421 of references[3].\\

Now, we give the theorem of this paper.\\

{\bfseries Theorem.}  If the function $ f(x) $  is defined in the interval $ [0,1]   $ and   satisfies the Dirichlet  condition, then,  when $ \chi_q(-1)=1  $,  we have

$$  \sum_{k=1}^{q-1 } \chi_q (k) f^* \left(\frac{k}{q}\right) = 2\tau (\chi_q) \sum_{n=1}^{\infty} \overline{\chi}_q (n) \int_0^1 f(t)\cos (2\pi nt) dt $$\\

when $ \chi_q(-1)=-1 $, we have\\

$$  \sum_{k=1}^{q-1 } \chi_q (k) f^* \left(\frac{k}{q}\right) =- 2i\tau (\chi_q) \sum_{n=1}^{\infty} \overline{\chi}_q (n) \int_0^1 f(t)\\sin (2\pi nt) dt $$\\

{\bfseries  Proof.}  By lemma 3,  we have

$$  f^*(x)=\int_0^1 f(t) dt +  \sum_{  n =-\infty \atop n  \neq 0 }^{\infty} \exp(2\pi i nx ) \int_0^1 f(t) \exp(- 2\pi i nt )  dt $$ \\

 we take $ x=\frac{k}{q} , \,  1 \leq  k\leq q-1 $,  then\\

$$  f^* \left(\frac{k}{q} \right)=\int_0^1 f(t) dt +  \sum_{  n =-\infty \atop n  \neq 0 }^{\infty} \exp \left(2\pi i \frac{ n k}{q} \right) \int_0^1 f(t) \exp(- 2\pi i nt )  dt $$ \\

Multiply the above formula by $ \chi_q (k)$, then sum over $ k $, we have\\

$$  \sum_{k=1}^{q-1} \chi_q (k)  f^* \left(\frac{k}{q} \right)=\sum_{k=1}^{q-1} \chi_q (k)  \int_0^1 f(t) dt  $$ \\

$$+  \sum_{  n =-\infty \atop n  \neq 0 }^{\infty} \sum_{k=1}^{q-1} \chi_q (k)   \exp \left(2\pi i \frac{ n k}{q} \right) \int_0^1 f(t) \exp(- 2\pi i nt )  dt $$ \\

By Lemma 1,

$$  \sum_{k=1}^{q-1} \chi_q (k)  \exp \left(2\pi i \frac{ n k}{q} \right)    = \overline{\chi}_q (n) \tau (\chi_q)   \quad  and    \quad  \sum_{k=1}^{q-1} \chi_q (k) = 0   $$  \\

We have

$$  \sum_{k=1}^{q-1} \chi_q (k)  f^* \left(\frac{k}{q} \right)=  \tau(\chi_q)  \sum_{  n =-\infty \atop n  \neq 0 }^{\infty} \overline{\chi}_q (n)  \int_0^1 f(t) \exp(- 2\pi i nt )  dt   $$ \\

$$  =  \tau(\chi_q)  \sum_{  n =1}^{\infty} \overline{\chi}_q (n)  \int_0^1 f(t) \exp(- 2\pi i nt )  dt  $$ \\

$$ +  \tau(\chi_q)  \sum_{ n =1}^{\infty} \overline{\chi}_q (-n)  \int_0^1 f(t) \exp( 2\pi i nt )  dt  $$\\

$$  =  \tau(\chi_q)  \sum_{  n =1}^{\infty} \overline{\chi}_q (n) \bigg ( \int_0^1 f(t) \exp(- 2\pi i nt )  dt   $$ \\

$$  +  \overline{\chi}_q (-1)  \int_0^1 f(t) \exp( 2\pi i nt  )  dt \bigg )  $$\\

therefore, when  $ \chi_q (-1)=1 $,  we have\\

$$ \sum_{k=1}^{q-1} \chi_q (k)  f^* \left(\frac{k}{q} \right)= 2 \tau(\chi_q)  \sum_{  n =1}^{\infty} \overline{\chi}_q (n)  \int_0^1 f(t) \cos( 2\pi  nt )  dt   $$ \\

when  $ \chi_q (-1)=-1 $,  we have\\

$$ \sum_{k=1}^{q-1} \chi_q (k)  f^* \left(\frac{k}{q} \right)=- 2i \tau(\chi_q)  \sum_{  n =1}^{\infty} \overline{\chi}_q (n)  \int_0^1 f(t) \sin ( 2\pi  nt )  dt   $$ \\

This completes the proof of the theorem.\\

From this theorem, we can see that the calculation of the character sum becomes the calculation of integrals.\\

Below, we give a few special examples.\\

{\bfseries  The first example.} \\

Let $ \chi_q $  be the primitive real character and  $ \chi_q (-1)=-1   $, then\\

$$ \sum_{k=1}^{q-1} \chi_q (k) \left(\frac{k}{q}\right)^2 = -\frac{\sqrt{q}}{\pi}  L(1,\chi_q)$$\\

By Theorem and Lemma 2,  easily seen\\

$$\sum_{k=1}^{q-1} \chi_q (k) \left(\frac{k}{q}\right)^2 \, = 2\sqrt{q}  \sum_{n=1}^\infty  \chi_q (n) \int_0^1  t^2 \, \sin (2\pi nt) dt $$\\

We compute the integral as follows\\

$$  \int_0^1  t^2 \, \sin (2\pi nt) dt = -\frac{1}{2\pi n}  \int_0^1  t^2 \, d \cos (2\pi nt)  $$

$$ =  -\frac{1}{2\pi n}  +  \frac{2}{2\pi n}  \int_0^1  t \, \cos (2\pi nt) dt =   -\frac{1}{2\pi n}  +  \frac{2}{(2\pi n)^2}  \int_0^1  t \, d \sin (2\pi nt)  $$\\

$$  =   -\frac{1}{2\pi n} -  \frac{2}{(2\pi n)^2}  \int_0^1 \, \sin (2\pi nt)  dt  =   -\frac{1}{2\pi n}  $$\\

therefore\\

$$ \sum_{k=1}^{q-1} \chi_q (k) \left(\frac{k}{q}\right)^2 =   -\frac{\sqrt{q}}{\pi} \sum_{n=1}^\infty \frac{\chi_q (n)}{n}   =  -\frac{\sqrt{q}}{\pi}  L(1,\chi_q)$$\\

{\bfseries The second example.} \\

Let  $ \chi_q $  be the primitive real character and $ \chi_q (-1)=1 $, then

$$ \sum_{k=1}^{q-1} \chi_q (k) \log k = - \frac{\sqrt{q}}{2} L(1,\chi_q) + c \sqrt{q} $$ \\

where $ c $ is a absolute constant.\\

By Theorem and Lemma 2,  we have\\

$$  \sum_{k=1}^{q-1} \chi_q (k) \log \frac{k}{q} = 2\sqrt{q} \sum_{n=1}^\infty \chi_q (n)  \int_0^1 \log t \cos (2\pi n t) dt  $$\\

Now,  let's compute the integral \\

$$    \int_0^1 \log t \cos (2\pi n t) dt =  \frac{1}{2\pi n }  \int_0^1 \log t \, d \sin (2\pi n t)  $$ \\

$$ = - \frac{1}{2\pi n}  \int_0^1 \frac{1}{t}  \sin (2\pi n t) dt = - \frac{1}{2\pi n}  \int_0^n \frac{1}{t}  \sin (2\pi  t) dt  $$\\

$$  =  - \frac{1}{2\pi n}  \int_0^\infty \frac{1}{t}  \sin (2\pi  t) dt +  \frac{1}{2\pi n}  \int_n^\infty \frac{1}{t}  \sin (2\pi  t) dt   $$\\

easily seen\\

$$ - \frac{1}{2\pi n}  \int_0^\infty \frac{1}{t}  \sin (2\pi  t) dt   =  - \frac{1}{2\pi n}  \int_0^\infty \frac{1}{t}  \sin t dt =  - \frac{1}{4 n} $$\\

because\\

 $$ \frac{1}{2\pi n}  \int_n^\infty \frac{1}{t} \,  \sin (2\pi  t) dt =  -  \frac{1}{(2\pi )^2 n }  \int_n^\infty \frac{1}{t } \, d \cos (2\pi  t)   $$\\

$$  =   \frac{1}{(2\pi n )^2 } -  \frac{1}{(2\pi )^2 n }  \int_n^\infty \frac{1}{t^2 } \,  \cos (2\pi  t) dt  \ll  \frac{1}{n^2} + \frac{1}{n}  \int_n^\infty \frac{1}{t^2 }  dt \ll  \frac{1}{n^2}   $$ \\

as well as\\

$$ \sum_{k=1}^{q-1} \chi_q (k) \log \frac{k}{q} = \sum_{k=1}^{q-1} \chi_q (k) \log k - \log q  \sum_{k=1}^{q-1} \chi_q (k) =   \sum_{k=1}^{q-1} \chi_q (k) \log k $$ \\

This completes the proof.\\

{\bfseries  The third example. }\\

Let $ \chi_q $  be the primitive real character and  $ \chi_q (-1)=1  $, then \\

$$  \sum_{k=1}^{q-1} \chi_q (k) \, e^{\frac{k}{q}} = 2(e-1) \sqrt{q}  \sum_{n=1}^\infty  \frac{\chi_q (n)}{ 1+ 4 \pi^2  n^2 } $$\\

Let $ \chi_q $  be the primitive real character and  $ \chi_q (-1)=-1  $, then \\

$$  \sum_{k=1}^{q-1} \chi_q (k) \, e^{\frac{k}{q}} =-4\pi(e-1) \sqrt{q}  \sum_{n=1}^\infty  \frac{\chi_q (n) \, n }{ 1+ 4 \pi^2  n^2 } $$\\

{\bfseries  Proof.}  When $ \chi_q (-1)=1  $, by Theorem and Lemma 2,  we have\\

$$ \sum_{k=1}^{q-1} \chi_q (k) \, e^{\frac{k}{q}} = 2\sqrt{q} \sum_{n=1}^\infty  \chi_q (n) \int_0^1 e^t \,  \cos (2\pi nt)  dt $$\\

Because, we know \\

$$ \int e^{ax} \cos bx \, dx= \frac{e^{ax}}{a^2+b^2} (a\cos bx+b\sin bx)  $$\\

Therefore\\

 $$  \int_0^1 e^{t} \cos (2\pi nt) \, dt=  \frac{e-1}{1+4\pi^2 n^2}   $$

Therefore\\

$$  \sum_{k=1}^{q-1} \chi_q (k) \, e^{\frac{k}{q}} =  2(e-1) \sqrt{q}  \sum_{n=1}^\infty  \frac{\chi_q (n)}{ 1+ 4 \pi^2  n^2 } $$\\

When $ \chi_q (-1)=-1  $, by Theorem and Lemma 2,  we have\\

$$ \sum_{k=1}^{q-1} \chi_q (k) \, e^{\frac{k}{q}} = 2\sqrt{q} \sum_{n=1}^\infty  \chi_q (n) \int_0^1 e^t \, \sin (2\pi nt)  dt $$\\

Because, we know \\

$$ \int e^{ax} \sin bx \, dx = \frac{e^{ax}}{a^2+b^2} ( a \sin bx - b \cos bx)  $$\\

Therefore\\

$$\int_0^1  e^{t} \sin (2\pi nt) \, dt = - (e-1) \frac{2\pi n}{1+4 \pi^2 n ^2}   $$\\

Therefore\\

$$  \sum_{k=1}^{q-1} \chi_q (k) \, e^{\frac{k}{q}} =-4\pi(e-1) \sqrt{q}  \sum_{n=1}^\infty  \frac{\chi_q (n) \, n }{ 1+ 4 \pi^2  n^2 } $$\\

This completes the proof.\\

{\bfseries  The fourth example. }\\

This is a well-known formula.\\

We write

$$ F(y) = \sum_{ 1 \leq  k \leq qy } \chi_q (k)  \qquad   and   \qquad  F^* (y) =  \frac{F(y+0)+F(y-0)}{2}   $$\\

When $ \chi_q (-1)=1 $,  we have\\

$$  F^* (y) = \frac{\tau (\chi_q)}{\pi}  \sum_{n=1}^\infty   \frac{ \overline{\chi}_q (n) }{n} \sin (2\pi n y) $$\\

When $ \chi_q (-1)=-1 $,  we have\\

$$  F^* (y) = \frac{\tau (\chi_q)}{ i \pi} L(1,\overline{\chi}_q)     - \frac{\tau (\chi_q)}{ i \pi}  \sum_{n=1}^\infty   \frac{ \overline{\chi}_q (n) }{n} \cos (2\pi n y) $$\\

{\bfseries  Proof. }   Let $ 0<y<1 $,  we define the function $ f(x) $  as follow\\

\[   f(x)  = \left\{
  \begin{array}{ll}
  1  \quad if \,\,\, 0\leq x \leq y \\
  0  \quad if \,\, \, y<x<1    \qquad \qquad \qquad
  \end{array}
  \right.
\]\\

By Theorem, when $ \chi_q (-1) =1 $,  we have\\

$$  \sum_{k=1}^{q-1} \chi_q (k) f^* \left(\frac{k}{q}\right ) = 2 \tau (\chi_q) \sum_{n=1}^\infty \overline{\chi}_q (n) \int_0^y  \cos (2\pi n t ) dt $$\\

Because\\

$$ \int_0^y  \cos (2\pi n t ) dt = \frac{1}{2\pi n} \int_0^y d \sin (2\pi n t ) =  \frac{ \sin (2\pi n y )  }{2\pi n}  $$\\

Therefore\\

$$  F^* (y) = \frac{\tau (\chi_q)}{\pi}  \sum_{n=1}^\infty   \frac{ \overline{\chi}_q (n) }{n} \sin (2\pi n y) $$\\

By Theorem,  when  $  \chi_q (-1) =-1 $,  we have\\

$$  \sum_{k=1}^{q-1} \chi_q (k) f^* \left(\frac{k}{q}\right ) =- 2 i \tau (\chi_q) \sum_{n=1}^\infty \overline{\chi}_q (n) \int_0^y  \sin (2\pi n t ) dt $$\\

Because\\

$$  \int_0^y  \sin (2\pi n t ) dt = - \frac{1}{2\pi n } \int_0^y d \cos (2\pi n t ) = - \frac{1}{2\pi n } ( \cos (2\pi n y ) -1 )  $$\\

Therefore \\

$$  F^* (y) = \frac{\tau (\chi_q)}{ i \pi} L(1,\overline{\chi}_q)     - \frac{\tau (\chi_q)}{ i \pi}  \sum_{n=1}^\infty   \frac{ \overline{\chi}_q (n) }{n} \cos (2\pi n y) $$\\

This completes the proof. \\

\vspace{10mm} \centerline{ REFERENCES } \vspace{5mm}

[1]  Hugh L. Montgomery, Robert C. Vaughan, { \itshape Multiplicative Number Theory I. Classical Theory,} Cambridge University Press, 2006. \\

[2]  Hua Loo Keng, {  \itshape  Introduction to Number Theory,  }  Springer-Verlag Berlin Heidelberg New York, 1982. \\

[3] I.N. Bronshtein, K.A.Semendyayev, G.Musiol, H.Muehlig,  {  \itshape   Handbook of Mathematics,  }  Springer Berlin Heidelberg New York 2005.

\end{document}